\documentclass{amsart}

\usepackage{amsmath}
\usepackage{amscd}
\usepackage{amssymb} 

\newcommand{\bb}{\mathbb}

\newcommand{\bC}{{\bb C}}

\newcommand{\bR}{{\bb R}}
\newcommand{\bZ}{{\bb Z}}
\newcommand{\cV}{{\mathcal V}}

\newcommand{\cVA}{{\mathcal V}{\mathcal A}}

\newcommand{\vac}{|0\rangle}

\DeclareMathOperator{\id}{id}

\DeclareMathOperator{\Aut}{Aut}

\DeclareMathOperator{\End}{End}
\DeclareMathOperator{\Ker}{Ker}
\DeclareMathOperator{\Img}{Im}

\newtheorem{theorem}{Theorem}[section]
\newtheorem{theorem/definition}{Theorem/Definition}[section]

\newtheorem{proposition}{Proposition}[section]
\newtheorem{lemma}{Lemma}[section]

\newtheorem{corollary}{Corollary}[section]

\theoremstyle{remark}
\newtheorem{remark}{Remark}[section]

\theoremstyle{definition}
 \newtheorem{example}{Example}[section]
\newtheorem{definition}{Definition}[section]

\begin{document}
\title{Superconformal vertex algebras in differential geometry. I.}
\author{Jian Zhou}
\address{Department of Mathematics\\
Texas A\&M University\\
College Station, TX 77843}
\email{zhou@math.tamu.edu}
\begin{abstract} 
We show how to construct an $N=1$ superconformal vertex algebra (SCVA)
from any Riemannian manifold.
When the Riemannian manifold has special holonomy groups,
we discuss the extended supersymmetry.
When the manifold is complex or K\"{a}hler,
we also generalize the construction to obtain $N=2$ SCVA's.
We study the BRST cohomology groups of the topological vertex algebras
obtained by the $A$ twist and the $B$ twist 
from these $N=2$ SCVA's. 
We show that for one of them,
the BRST cohomologies are isomorphic to 
$H^*(M, \Lambda^*(T^*M))$ and $H^*(M, \Lambda^*(TM))$ respectively.
This provides a mathematical formulation of the $A$ theory and $B$ theory 
in physics literature.
The connection with elliptic genera is also discussed.
Furthermore,
when the manifold is hyperk\"{a}hler,
we generalize our constructions to obtain $N=4$ SCVA's.
A heuristic relationship with super loop space is also discussed.
\end{abstract}
\date{\today}
\maketitle

\section{Introduction}

The goal of this paper is to associate superconformal vertex algebras
to smooth manifolds, 
especially when the manifolds are equipped with various structures 
as indicated in the abstract.
Such results are important in the author's program of studying mirror symmetry
in terms of Quillen's closed model category theory \cite{Qui}.

It has been known for some time in the physics literature 
that the sigma model has an $N=1$ supersymmetric extension 
when the target manifold is Riemannian,
an $N=2$ supersymmetric extension if the target manifold
is K\"{a}hler,
and an $N=4$ supersymmetric extension if the target manifold
is hyperk\"{a}hler.
See  e.g. Zumino \cite{Zum}, Alvarez-Gaum\'e and Freedman \cite{Alv-Fre}.
In the Riemannian case,
when the holonomy group of the target manifold is $G_2$ or $Spin(7)$,
the supersymmetry can be extended 
(see e.g. Shatashvilli and Vafa \cite{Sha-Vaf});
in the K\"{a}hler case,
when the target manifold $M$ has holonomy group $SU(n)$ 
(i.e. $M$ is Calabi-Yau),
the extended supersymmetry was considered earlier by Odake \cite{Oda}.
Recall that the sigma model involves the maps from 
a two-dimensional Riemannian manifold to 
another Riemannian manifold $M$ (the target manifold).
Since the type of supersymmetric extension depends on the target manifold,
it is reasonable that the target manifold alone is the source of the
supersymmetry.
If that is the case, 
one should be able to obtain supersymmetry just from the target manifold,
and deduce from it the supersymmetry of the sigma model. 
In this paper,
we construct various bundles of superconformal vertex algebras
from smooth manifolds,
depending on whether they admit
Riemannian metrics with various holonomy groups,
or complex structures.
We show that the spaces of sections of these bundles
are superconformal vertex algebras of the same types as their fibers.
In the cases of complex, K\"{a}hler, and hyperk\"{a}hler manifolds,
we also show that these bundles can be regarded as holomorphic vector bundles
hence one can consider the $\bar{\partial}$ operators on them.
The corresponding Dolbeault cohomologies are shown to be 
superconformal vertex algebras,
again of the same types as the fibers of the relevant bundles.
Our constructions should lead to rigorous explanation of the
supersymmetric extensions of sigma models with Riemannian manifolds 
(having special holonomy groups) as target manifolds.
Notice that sigma models are not explicitly involved in \cite{Sha-Vaf}.
Our constructions might be implicitly known to the physicists.
But to the author's best knowledge,
they have not been explicitly formulated in the present form.

Part of the motivations of this paper is to provide mathematical 
formulations of the $A$ theory and the $B$ theory 
without using the sigma model.
A very important notion in string theory is that of an $N=2$ superconformal 
field theory (SCFT).
We understand such a theory as 
the assignment of $N=2$ SCVA's to suitable physical situations.
Physicists showed that the primary chiral states of an
$N=2$ SCVA form an algebra.
Similarly for the primary anti-chiral states.
See e.g. Lerche-Vafa-Warner \cite{Ler-Vaf-War}.
See also Warner \cite{War} for an exposition.
A closely related notion is that of a topological vertex algebra
of which one can consider the BRST cohomologies.
Given an $N=2$ SCVA,
there are two ways to twist it 
(called the $A$ twist and the $B$ twist respectively)
to obtain a topological vertex algebra
(Eguchi and Yang \cite{Egu-Yan}).
The BRST cohomology groups of these two topological vertex algebras
correspond to the algebras of the primary chiral and anti-chiral states 
of the original $N=2$ SCVA respectively.
Given a Calabi-Yau manifold $M$,
it has been widely discussed in physics literature for many years
that there is an $N=2$ SCFT associated to it,
with the two twists giving the $A$ theory and the $B$ theory 
respectively.
Witten \cite{Wit6} discuss such a theory starting with sigma models with
a Calabi-Yau manifold as the target manifold.
However,
he did not explicitly give the construction of an $N=2$ SCVA.
The primary chiral algebra 
(i.e. the BRST cohomology in the $A$ theory)
is argued to be isomorphic to the cohomology of the manifold
as vector spaces,
but not as algebras,
hence leading to the notion of quantum cohomology 
(\cite{Ler-Vaf-War}, \cite{Wit6}).
Since the fundamental work by Ruan and Tian \cite{Rua-Tia},
the mathematical theory of quantum cohomology motivated by the sigma model  
has reached maturity through
the efforts of mathematicians too many to list here.
However,
since such developments do not require the $N=2$ superconformal
vertex algebras in the discussions,
the problem of associating $N=2$ SCVA to Calabi-Yau manifolds
has not been fully addressed in the mathematical literature.
We will fill this gap in this paper,
using constructions not involving sigma models. 

It is a good time to compare with earlier work by  
Malikov, Schechtman, and Vaintrob \cite{Mal-Sch-Vai}.
They constructed 
for any Calabi-Yau manifold a sheaf of topological vertex algebras.
Their theory corresponds to the $A$ theory.
It is not clear how to extend their construction to encompass the $B$ theory.
The first essential difference between their theory and ours is that
we construct an $N=2$ SCVA from any complex manifold $M$,
whose twists give $H^{*,*}(M)$ and $H^{-*,*}(M)$ 
when taking the BRST cohomologies.
The $A$ twist and the $B$ twist correspond 
to the $A$ theory and the $B$ theory in the case of Calabi-Yau manifolds,
but our construction is more general.
As mentioned above,
we give a unified treatment for Riemannian manifolds 
with holonomy groups of all possible kinds,
and for complex manifolds.
The second difference is that we use bundles of vertex algebras while
they use sheaves of vertex algebras
that do not come from vector bundles.
Our constructions have the advantage of being easy and admitting
treatments by usual techniques in differential geometry,
while their construction involves complicated localization argument.

The connection with elliptic genera is manifest in our approach.
As a byproduct,
we introduce the $K$-theory of vertex algebra bundles 
and its conformal and superconformal cousins.
Such theories seem to deserve further investigation 
because of possible links with the theory of elliptic cohomology.
Note that the $K$-theories of infinite dimensional vector bundles
was first introduced by Liu \cite{Liu} in the search of a geometric
description of elliptic cohomology.
The bundles of vertex algebras have been used by Tamanoi \cite{Tam}.

We remark that we have only used the construction of vertex algebras
from a vector space with an inner product in this paper.
Other constructions of vertex algebras can be applied in the same fashion.
In particular,
we will apply the constructions of vertex algebras
from a symplectic vector space to symplectic geometry
is a sequel \cite{Zho}.
For example,
we will associate $N=2$ SCVA's to a symplectic manifold with a Lagrangian
fibration,
hence establish some relationship between our work and 
the Strominger-Yau-Zaslow program \cite{Str-Yau-Zas}.

The rest of the paper is arranged as follows.
In Section 2 we recall some basic definitions and facts about vertex
algebras, including
the definitions of $N=2$ superconformal 
vertex algebras and their primary (anti-)chiral algebras,
and their relationship with topological vertex algebras and BRST cohomology.
We define vertex algebra bundles and their $K$-theories in Section 3. 
In Sectio 4 we apply the construction of an $N=1$ SCVA from a vector space 
with an inner product to the tengent bundle of Riemnnian manifold.
We explain extended supersymmetry discovered by physicists
when the Riemanian manifold has special holonomy group
in our framework.
Section 5 applies the construction of an $N=2$ SCVA from a vector space 
with an inner product and a polarization to the tangent bundle 
of a complex or K\"{a}hler manifold.
Our main result is Theorem \ref{thm:N=2complex} 
and Theorem \ref{thm:N=2Kahler}.
We apply in Section 6 the construction of an $N=4$ SCVA from a quternionic
vector space to the tangentbundle of a hyperk\"{a}hler manifold.
Finally,
we explain our choice of the vector bundles in terms of loop spaces 
in Section 7.

\section{Preliminaries on vertex algebras}

\subsection{Definition of a conformal vertex algebra}

A {\em field} acting on a vector space $V$ is a formal power series
$$a(z) = \sum_{n \in \bZ} a_{(n)}z^{-n-1} \in \End(V)[[z, z^{-1}]]$$ 
such that such that for any $v \in V$,
$a_{(-n)}v = 0$ for $n \gg 0$.
A {\em vertex algebra} consists of the following data:
\begin{itemize}
\item[(i)] a $\bZ$-graded vector space $V$ (the space of states),
\item[(ii)] a vector $\vac \in V_0$ (the vacuum vector),
\item[(iii)] a degree zero linear map $Y: V \to \End V[[z, z^{-1}]]$ 
(the state-field correspondence) whose image lies in the space of fields,
denoted by $Y(a, z) = \sum_{n \in \bZ} a_{(n)}z^{-n-1}$,
\end{itemize} 
satisfying the following axioms:
\begin{itemize}
\item[(1)] (translation covariance): $[T, Y(a, z)] = \partial_zY(a, z)$,
where $T \in \End V$ is defined by $T(a) = a_{-2}\vac$ for $a \in V$,
\item[(2)] (vacuum): $Y(\vac, z) = \id_V$, 
$Y(a, z)\vac|_{z = 0} = a$,
\item[(3)] (locality):
for any $a, b \in V$,
$Y(a, z)$ and $Y(b, w)$ are mutually local.
\end{itemize}

\begin{example}
A vertex algebra $V$ is called {\em holomorphic} 
if $a_{(n)} = 0$ for $n \geq 0$.
Holomorphic vertex algebras correspond to graded commutative associative 
algebras with unit and with a derivation $T$ of even degree:
$$Y(a, z)b = e^{zT}(a) \cdot b.$$
\end{example}

A {\em conformal vertex algebra of rank $c$} is a vertex algebra $V$ with
a {\em conformal vector} $\nu$ of degree $0$
such that the corresponding field $Y(\nu, z) = \sum_{n \in \bZ} L_n z^{-n-2}$
satisfies:
\begin{eqnarray} \label{eqn:Virasoro}
[L_m, L_n] = (m-n) L_{m+n} + \frac{(m^3-m)c}{12}\delta_{m+n, 0};	
\end{eqnarray}
furthermore, $L_{-1} = T$,
and $L_0$ is diagonizable.
Some authors have used the term {\em vertex operator algebras}.

\subsection{$N=1$ superconformal vertex algebras}

An {\em $N = 1$ superconformal vertex algebra}
a vertex algebra $V$ of order $c$ with an
odd vector $\tau$ (called {\em $N=1$ superconformal vector}), 
such that the field $G(z) = Y(\tau, z) 
= \sum_{n \in \frac{1}{2} + \bZ} G_m z^{-n -\frac{3}{2}}$ satisfies 
\begin{eqnarray*}
&& [G_m, G_n] = 2L_{m+n} + \frac{1}{3}(m^2- \frac{1}{4})\delta_{m+n}c, \\
&& [G_m, L_n] = (m- \frac{n}{2})G_{m+n}
\end{eqnarray*}

\subsection{$N =2$ superconformal vertex algebras}

An  {\em $N=2$ superconformal vertex algebra} (SCVA)
is a conformal vertex algebra $V$ with two odd vectors $\tau^{\pm}$ 
and an even vector $j$ such that the fields
$G^{\pm}(z) = Y(\tau^{\pm}, z)
= \sum_{r \in \frac{1}{2}+\bZ} G^{\pm}_r z^{-r-\frac{3}{2}}$ 
and $J(z) = Y(j, z) = \sum_{n \in \bZ} J_nz^{-n-1}$ satisfy 
\begin{eqnarray}
& [G_r^-, G_s^+] = 2L_{r+s} - (r -s)J_{r+s} 
	+ \frac{c}{3}(r^2 - \frac{1}{4}) \delta_{r+s, 0}, 
	\label{eqn:GG} \\
& [G_r^+, G_s^+] = [G_r^-, G_s^-] = 0, \label{eqn:GG2}\\
& [L_m, L_n] = (m-n)L_{m+n} + \frac{c}{12}m(m^2-1)\delta_{m+n, 0}, 
	\label{eqn:LL} \\
& [L_n, G_r^{\pm}] = (\frac{n}{2} -r)G_{n+r}^{\pm}, \label{eqn:LG}\\
& [L_n, J_m] = - mJ_{m+n}, \label{eqn:LJ}\\
& [J_m, J_n] = \frac{c}{3}m\delta_{m+n, 0}, \label{eqn:JJ}\\
& [J_n, G_r^{\pm}] = \pm G^{\pm}_{n+r}.\label{eqn:JG}
\end{eqnarray}

It is easy to see that for any nonzero $a \in \bC$,
$a\tau^+ + \frac{1}{a}\tau^-$ is an $N=1$ superconformal vector.

\subsection{Primary chiral algebra of an $N=2$ SCVA}

A state $a$ in an $N=2$ superconformal vertex algebra is called
{\em primary} of conformal weight $h$ and $U(1)$ charge $q$ 
if 
\begin{eqnarray*}
&& L_na = J_n a = 0, \;\;\;\; n \geq 1, \\
&& L_0a = ha,  \;\;\;\; J_0 a =qa, \\
&& G_r^{\pm}a = 0, \;\;\;\; r \geq \frac{1}{2}. 
\end{eqnarray*}

An $N=2$ SCVA $V$ is said to be {\em unitary} if there a positive definite
Hermitian metric
$\langle \cdot|\cdot\rangle$ on $V$ such that $(G_r^+)^* = G_{-r}^-$. 
It is easy to see that 
if $(V, \langle\cdot|\cdot\rangle)$ is a unitary $N=2$ SCVA,
then we have
\begin{align*}
L_n^* & = L_{-n}, & J_n^* = J_{-n}.
\end{align*}
An $N=2$ SCVA is said to be {\em nondegenerate} 
if the eignevalues of $L_0$ are discrete, bounded from below, and with finite
dimensional eigenspaces.

In an $N=2$ SCVA, 
a state $a$ is called {\em chiral} if it satisfies:
$$G^+_{-\frac{1}{2}}a = 0.$$
Similarly,
a field $a(z)$ is called {\em chiral} if it satisfies:
$$[G^+_{-\frac{1}{2}}, a(z)] = 0.$$
Anti-chiral states and fields are defined with $+$ replaced by $-$.
We now recall some important results from 
Lerche-Vafa-Warner \cite{Ler-Vaf-War}.

\begin{lemma} \label{lm:weightcharge}
In a unitary $N=2$ SCVA $V$,
if $a$ is a vector of conformal weight $h$ and $U(1)$ charge $q$,
then 
$$h \geq |q|/2,$$
with $h = q/2$ (resp. $h = - q/2$) 
iff $a$ is a primary chiral (resp. anti-chiral) state.
\end{lemma}

\begin{corollary}
In a unitary $N=2$ SCVA $V$,
the primary chiral states form an graded commutative associative algebra 
induced by the normally ordered product:
the product between states $a$ and $b$ is given by
 $a_{(-1)}b$.
Similarly for the primary anti-chiral states.
(These algebras will be referred to as the {\em primary chiral algebra} 
and {\em primary anti-chiral algebra} of $V$ respectively.)
\end{corollary}

\begin{lemma}
Let $V$ be a unitary $N=2$ SCVA of central charge $c$.
if $a \in V$ is of conformal weight $h$ and $U(1)$ charge $q$,
then one has
$$h \leq c/6.$$
\end{lemma}

\begin{corollary}
The  space of primary chiral states 
in a nondegenerate unitary $N=2$ SCVA is finite dimensional.
\end{corollary}

\subsection{Topological vertex algebras and BRST cohomology}

Closely related to $N=2$ SCVA are the {\em topological vertex algebras}.
Recall that a topological vertex algebra of rank $d$ 
is a conformal vertex algebra of central charge $0$,
equipped with an even element $J$ of conformal weight $1$,
an odd element $Q$ of conformal weight $1$,
and an odd element $G$ of conformal weight $2$,
such that their fields satisfy the following OPE's:
\begin{eqnarray}
&& T(z)T(w) \sim \frac{2T(w)}{(z-w)^2} + \frac{\partial_W T(w)}{z-w}, 
	\label{OPE:TTT}\\
&& J(z)J(w) \sim \frac{d}{(z-w)^2}, \label{OPE:TJJ} \\
&& T(z)J(w) \sim - \frac{d}{(z-w)^3} + \frac{J(w)}{(z-w)^2} 
	+ \frac{\partial_w J(w)}{z-w},  \label{OPE:TTJ} \\
&& G(z)G(w) \sim 0, \label{OPE:TGG}\\
&& T(z)G(w) \sim \frac{2G(w)}{(z-w)^2} + \frac{\partial_w G(w)}{z-w},
	\label{OPE:TTG}\\
&& J(z)G(w) \sim - \frac{G(w)}{z-w},  \label{OPE:TJG}\\
&& Q(z)Q(w) \sim 0, \label{OPE:TQQ} \\
&& T(z)Q(w) \sim \frac{Q(w)}{(z-w)^2} + \frac{\partial_w Q(w)}{z-w},  
	\label{OPE:TTQ}\\
&& J(z)Q(w) \sim \frac{Q(w)}{z-w},  \label{OPE:TJQ}\\
&& Q(z)G(w) \sim \frac{d}{(z-w)^3} + \frac{J(w)}{(z-w)^2} + \frac{T(w)}{z-w}. 
	\label{OPE:TQG}
\end{eqnarray}
Here we have written the Virasoro field associated with the conformal vector
as $T(z) = \sum_{n \in \bZ} T_n z^{-n-2}$.
As a consequence of (\ref{OPE:TQQ}),
$$Q_0^2 = \frac{1}{2}[Q_0, Q_0] = 0.$$
The operator $Q_0$ is called the {\em BRST operator},
and the cohomology
$$H^*(V, Q_0) = \Ker Q_0/\Img Q_0$$
is called the {\em BRST cohomology}.
From (\ref{OPE:TQQ}),
one gets
$$T(z) = [Q_0, G(z)].$$
Hence for any $v \in \Ker Q_0$,
we have
$$T_n v = [Q_0, G_n]v = Q_0G_nv \in \Img Q_0.$$
In other words,
the map $[v] \mapsto [T_n v]$ induces a trivial representation of the Virasoro
algebra on the BRST cohomology.

Inspired by Witten \cite{Wit1, Wit2},
Eguchi and Yang \cite{Egu-Yan} discovered the following 
important twisting construction:

\begin{proposition} \label{prop:top}
Given an $N=2$ SCVA $V$ with Virasoro field $L(z)$,
supercurrents $G^{\pm}(z)$ and $U(1)$ current $J(z)$,
one obtains a topological vertex algebra by taking:
\begin{align*}
& T(z) = L(z) + \frac{1}{2} \partial_z J(z), &&  
J_{top}(z) = J(z), \\
& Q(z) = G^+(z), &&  
G(z) = G^-(z), 
\end{align*}
or
\begin{align*}
& T(z) = L(z) - \frac{1}{2} \partial_z J(z), && 
J_{top}(z) = -J(z), \\
& Q(z) = G^-(z), &&
G(z) = G^+(z).
\end{align*}
Conversely,
given a topological vertex algebra,
one can obtain an $N=2$ SCVA structure on it by
\begin{align*}
& L(z) = T(z) - \frac{1}{2} \partial_z J_{top}(z), && 
J(z) = J_{top}(z), \\
& G^+(z) = Q(z), && 
G^-(z) = G(z), 
\end{align*}
or 
\begin{align*}
& L(z) = T(z) + \frac{1}{2} \partial_z J_{top}(z), &&
J(z) = -J_{top}(z), \\
& G^+(z) = G(z), && 
G^-(z) = Q(z).
\end{align*}
In the above,
we have used $J_{top}$ to denote the $U(1)$ charge for the topological
vertex algebra.
\end{proposition}

\begin{definition}
The two twists in Proposition \ref{prop:top} will be referred to as 
the $A$ twist and the $B$ twist respectively.
\end{definition}

As remarked in Lian-Zuckerman \cite{Lia-Zuc}, \S 3.9.4,
the BRST cohomology of a topological vertex algebra is 
graded commutative and associative.
The following results from Lerche-Vafa-Warner \cite{Ler-Vaf-War}
provide an alternative explanation.

\begin{lemma} (Hodge decomposition)
In a unitary $N=2$ SCVA $V$ of central charge $c$,
any state of conformal weight $h$ and $U(1)$ charge $q$ 
can be uniquely written as
$$a = a_0 + G^+_{-\frac{1}{2}}a_+ + G^-_{\frac{1}{2}}a_-,$$
for some primary chiral state $a_0$ and some states $a_+$ and $a_-$.
Furthermore,
when $a$ is chiral,
then one can take $a_- = 0$.
\end{lemma}

\begin{proposition}
For a nondegenerate unitary $N=2$ SCVA,
the primary chiral (resp. anti-chiral) algebra 
is isomorphic to the BRST cohomology
of the $A$ twist (resp. the $B$ twist).
\end{proposition}

\section{Vertex algebra bundles and the related $K$-theories}

\subsection{The category of vertex algebras}

We recall some categorical notions which will be used later.
A {\em homomorphism} from a vertex algebra $V$ to a vertex algebra $V'$ 
is a linear map $f: V \to V'$ of degree $0$ such that 
$$f(a_{(n)}b) = f(a)_{(n)}f(b)$$
for all $a, b \in V$, $n \in \bZ$.
One can define a category $\cVA$ which has objects the vertex algebras 
and as morphisms the vertex algebra homomorphisms.
Isomorphisms and automorphisms in this category are defined as usual.
Given two vertex algebras $U$ and $V$,
$U \oplus V$ is a vertex algebra with
\begin{eqnarray*}
&& \vac_{U \oplus v} = \vac_U \oplus \vac_V, \\
&& Y(u \oplus v, z) = Y(u,z) \oplus Y(v, z) 
= \sum_{n\in\bZ} (u_{(n)} \oplus v_{(n)}) z^{-n-1}.
\end{eqnarray*}
The direct sum gives the (categorical) product in the category $\cVA$.
Similarly,
$U \otimes V$ is a vertex algebra with
\begin{eqnarray*}
&& \vac_{U \otimes V} = \vac_U \otimes \vac_V, \\
&& Y(u \otimes v, z) = Y(u,z) \otimes Y(v, z) 
= \sum_{m,n\in\bZ} u_{(m)} \otimes v_{(n)} z^{-m-n-2}.
\end{eqnarray*}
This tensor product and the natural isomorphism $U \otimes V \to V \otimes U$
give $\cVA$ the structure of a symmetric monoidal category
(see e.g. Mac Lane \cite{Mac} for definition).

Other categorical notions such as subobjects and quotient objects 
can also be easily defined.
A {\em subalgebra} of a vertex algebra $V$ is a subspace $U$ containing $\vac$ 
such that $a_{(n)}U \subset U$
for all $a \in U$.
It is clear that $U$ is a vertex algebra with
$Y(a, z) = \sum_{n \in \bZ} a_{(n)}|_U z^{-n-1}$.
An {\em ideal} of a vertex algebra $V$ is a $T$-invariant subspace $I$ 
not containing $\vac$ such that $a_{(n)}I \subset I$ for all $a \in V$.
By skew symmetry,
one sees that
$$a_{(n)}V \subset I$$
for all $a \in I$.
Hence the quotient space $V/I$ has a structure of a vertex algebra such that
the quotient map is a homomorphism of vertex algebras.

The above definitions can be easily extended to the case of
charged vertex algebras, conformal vertex algebras 
and $N=n$ superconformal vertex algebras.
We remark on the case of tensor product.
Let $U$ be a vertex algebra and $V$ be a conformal vertex algebra 
with conformal vector $\nu_V$.
Then $U \otimes V$ is a conformal vertex algebra with conformal vector
$\vac_U \otimes \nu_V$.
This makes the category of conformal vertex algebras a {\em module}
over the category of vertex algebra.
Of course when $U$ is also conformal with conformal vector $\nu_U$,
then we can take $\nu_U \otimes \vac_V + \vac_U \otimes \nu_V$ to be
the conformal vector of $U \otimes V$.
Similarly for $N=n$ superconformal vertex algebras.

\subsection{Vertex algebras bundles}

Let $V$ be a vertex algebra,
denote by $\Aut(V)$ the automorphism group of $V$.
Let $M$ be a smooth topological space,
a vertex algebra bundle with fiber $V$ over $M$ is a vector bundle
$\pi: E \to M$ with fiber $V$ 
such that the transition functions lie in $\Aut(V)$.
Similarly define conformal vertex algebra bundles 
and superconformal vertex algebra bundles.
When $M$ is a smooth manifold or a complex manifold,
one can also define smooth or holomorphic vertex algebra bundles.

\begin{remark}
Given a vertex algebra $V$,
it is interesting to construct a universal vertex algebra bundle
$\cV \to B_V$ with fiber $V$ such that every other vertex algebra
bundle $E \to M$ with fiber $V$ can be obtained by pulling back
along a continuous map $f: M \to B_V$.
\end{remark}

\begin{lemma} \label{lm:section}
Given a vertex algebra bundle $E \to M$,
the space $E(M)$ of sections has an induced structure of a vertex algebra.
Similarly for 
(charged) conformal vertex algebra bundles and superconformal vertex 
algebra bundles.
\end{lemma}

\begin{proof}
Since the vacuum $\vac$ is preserved by the automorphisms,
it defines a section which we denote by $\vac_M$.
Given two sections $A$ and $B$,
the assignment
$$x \in M \mapsto A(x)_{(n)}B(x)$$
defines a section denoted by $A_{(n)}B$.
It is straightforward to check that 
$$Y(A,z)B = \sum_{n \in \bZ} A_{(n)}B z^{-n-1}$$
then defines a structure of a vertex algebra on $E(M)$.
\end{proof}

\subsection{Differential vertex algebras and their cohomology algebras}

A {\em derivation} of degree $k$ on a vertex algebra $V$ is a linear map
$S: V \to V$ of degree $k$ such that
$$S(a_{(n)}b) = (Sa)_{(n)}b + (-1)^{k|a|}a_{(n)}(Sb)$$
for all $a \in V$, $n \in \bZ$.
A {\em differential} is a derivation $d$ of degree $1$ such that $d^2 = 0$.
One can easily extend the above definitions to conformal 
vertex algebras and $N=n$ SCVA's.
For example,
in the case of conformal vertex algebras,
a derivation is required to annihilate the conformal vector $\nu$.

We omit the routine proof of the following:

\begin{lemma} \label{lm:cohomology}
Given a differential $d$ on a vertex algebra $V$,
$\Ker d$ is a vertex subalgebra of $V$,
$\Img d$ is an ideal of $\Ker d$.
Hence the cohomology $H(V, d) = \Ker d/\Img d$ is a vertex algebra.
We have similar results for (charged) conformal 
vertex algebras and $N=n$ SCVA's.
\end{lemma}

\subsection{$K$-theories of vertex algebra bundles}

The isomorphism classes of vertex algebra bundles 
on a topological space $X$ form an abelian monoid under the direct sum.
Denote by $K^{va}(X)$ the Grothendieck group of this monoid.
The tensor product induces a structure of a ring on this group.
We can similarly define $K$ theory of
conformal vertex algebra bundles, charged conformal vertex algebra bundles,
and $N = n$ superconformal vertex algebra bundles.
We denote them by $K^{cva}(X)$, $K^{ccva}(X)$, and $K^{n-scva}(X)$ respectively.

\begin{remark}
The $K$ theories of infinite dimensional vector bundles were introduced
by Liu \cite{Liu} in the search of geometric constructions of elliptic
cohomology.
The bundles of vertex algebras have been used by Tamanoi \cite{Tam}.
Our construction is partly inspired by their work.
\end{remark}

We restrict our attentions to conformal vertex algebra bundles
such that in each fiber $L_0$ is diagonalizable with finite dimensional 
eigenspaces.
For such bundles,
we consider
$$G_q(E) = \sum_n E_{c_n, h_n}q^{h_n - \frac{c_n}{24}},$$
where $E_{c_n, h_n}$ is the subbundle of elements of central charge $c_n$
and conformal weight $h_n$.
This extends to a homomorphism of the corresponding $K$-theories.
For an $N=2$ SCVA bundle $E$,
we also define
$$G_{q, y}(E) = \sum_n E_{c_n, h_n, j_n}q^{h_n - \frac{c_n}{24}}y^{j_n},$$
where $E_{c_n, h_n, j_n}$ is the subbundle of elements of central charge $c_n$,
conformal weight $h_n$, and $U(1)$ charge $j_n$.
This extends to a homomorphism of the corresponding $K$-theories.

\section{$N=1$ SCVA bundle from a Riemannian manifold}

\subsection{$N=1$ SCVA from a vector space with an inner product}
Let $T$ be a finite dimensional complex vector space with 
an inner product $g: T \otimes T \to \bC$.
The space
$$B(T, g) = S(\oplus_{n < 0}t^nT) = \otimes_{n > 0} S(t^{-n}T)$$
is spanned by elements of the form
$$a^1_{-j_1-1} \cdots a^m_{-j_m-1},$$
where $a^1, \cdots, a^m \in T$,
and for any $a \in T$, $n \in \bZ$,
$a_n = t^na$.
For any $a \in V$, 
define a field on $B(T, g)$ by
$$a(z) = \sum_{n \in \bZ} a_{(n)}z^{-n-1},$$
where for $n < 0$,
$a_{(n)}$ is symmetric product by $a_n$; 
for $n \geq 0$, $a_{(n)}$ is $n$ times the contraction by $a_{-n}$.
Similarly,
the space 
$$F(T, g)_{NS} = \Lambda(\oplus_{n > 0} t^{-n+\frac{1}{2}}T)
= \otimes_{n > 0} \Lambda(t^{-n+\frac{1}{2}}T)$$
is spanned by elements of the form
$$\phi^1_{-j_1-\frac{1}{2}} \cdots \phi^m_{-j_m-\frac{1}{2}},$$
where $\phi^1, \cdots, \phi^m \in T$,
$j_1, \cdots, j_m \geq 0$.
For $\phi \in T$,
$\phi_{-n+\frac{1}{2}}$ stands for $t^{-n+\frac{1}{2}}\phi$.
Define a field
$$\phi(z) = \sum_{r \in \frac{1}{2} + \bZ} \phi_{(r)} z^{-r - \frac{1}{2}},$$
where for $r < 0$,
$\phi_{(r)}$ is the exterior product by $\phi_r$;
for $r> 0$,
$\phi_{(r)}$ is the contraction by $\phi_{-r}$.

\begin{proposition} \label{prop:N=1}
(a) There is a structure of conformal vertex algebra on 
$B(T, g)$ defined by
$$Y(a^1_{-j_1-1} \cdots a^m_{-j_m-1}, z)
= :\partial^{(j_1)}a^1(z) \cdots \partial^{(j_m)}a^m(z):$$
for $a^1, \cdots, a^m \in T$ and $j_1, \cdots, j_m \geq 0$,
with conformal vector
$$\nu_B = \frac{1}{2}\sum_i a^i_{-1}a^i_{-1}$$
of central charge $c= \dim T$,
where $\{a^i\}$ is an orthonormal basis of $T$.

(b) There is a structure of a conformal vertex algebra on $F(T, g)_{NS}$ 
defined by
$$Y(\phi^1_{-j_1-\frac{1}{2}} \cdots \phi^m_{-j_m-\frac{1}{2}}, z)
= :\partial^{(j_1)}\phi^1(z) \cdots \partial^{(j_n)}\phi^m(z):$$
for $\phi^1, \cdots, \phi^m \in T$ 
and integers $j_1, \cdots, j_m \geq 0$,
with conformal vector
$$\nu_F = \frac{1}{2}\sum_i \phi^i_{-\frac{3}{2}}\phi^i_{-\frac{1}{2}}$$
of central charge $c = \frac{1}{2} \dim T$,
where $\{\phi^i\}$ is an orthonormal bases of $T$.

(c) The tensor product $V(T, g) = B(T, g) \otimes F(T, g)_{NS}$
is an $N = 1$ SCVA of central charge $\frac{3}{2}\dim T$ with 
\begin{align*}
\tau & = \sum_i a^i_{-1}\phi^i_{-\frac{1}{2}}, &
\nu & = \frac{1}{2} \sum_i a^i_{-1}a^i_{-1} 
+ \frac{1}{2} \sum_i \phi^i_{-\frac{3}{2}}\phi^i_{-\frac{1}{2}},
\end{align*}
where $\{e^i\}$ is an orthonormal basis of $T$,
and we write it as $\{a^i\}$ for the copy of $T$ in $B(T, g)$,
and as $\{\phi^i\}$ for the copy of $T$ in $F(T, g)_{NS}$.
\end{proposition}

\begin{remark} \label{rmk:N=1}
It is clear that $\nu_B$, $\nu_F$ and $\nu$ are independent of the choice of
the basis.
Furthermore, $O(T, g)$ acts on $B(T, g)$, $F(T, g)$ and $V(T, g)$ 
as automorphisms.
\end{remark}

\subsection{$N=1$ SCVA bundles from Riemannian manifolds}
\label{sec:Riemannian}

For any Riemannian manifold $(M, g)$,
we consider the principal bundle $O(M, g)$ of orthonormal frames.
Pick a point $x \in M$.
The structure group of $O(T, g)$ is $O(T_xM, g_x)$,
which acts on $V(T_xM \otimes \bC, g \otimes \bC)$ by automorphisms.
Applying Proposition \ref{prop:N=1},
Remark \ref{rmk:N=1} and Lemma \ref{lm:section},
we get the following:

\begin{theorem} \label{thm:N=1Riemannian}
Let $(M, g)$ be a Riemannian manifold.
Then 
$$V(TM\otimes \bC, g\otimes \bC)_{NS} 
= \otimes_{n>0} \Lambda(t^{-n+\frac{1}{2}}TM \otimes \bC) 
	\otimes_{n>0} S(t^{-n}TM\otimes \bC)$$
is an $N=1$ SCVA bundle,
hence $\Gamma(M, V(TM, g)_{NS})$ is an $N=1$ SCVA.
\end{theorem}

The bundle $V(TM\otimes \bC, g\otimes \bC)_{NS}$ has appeared in
the theory of elliptic genera.
It is easy to see that
$$G_q(V(TM\otimes \bC, g\otimes \bC)_{NS}
= q^{-\frac{\dim T}{16}} \otimes_{n > 0} S_{q^n} (TM \otimes \bC) 
\otimes_{n > 0} \Lambda_{q^{n -\frac{1}{2}}}(TM \otimes \bC)$$
(cf. Witten \cite{Wit2}, (27)).
A related formal power series
$$q^{-\frac{\dim T}{16}} \otimes_{n > 0} S_{q^n} (TM \otimes \bC) 
\otimes_{n > 0} \Lambda_{-q^{n -\frac{1}{2}}}(TM \otimes \bC)$$
(cf. Liu \cite{Liu})
can be obtained by introducing an operator $(-1)^F$.

\subsection{Riemannian manifolds 
of special holonomy groups and extended supersymmetries}

By a theorem of Berger \cite{Ber} (see also Simons \cite{Sim}),
the local holonomy groups of a nonsymmetric Riemannian manifold
can only be $O(n)$, $SO(n)$, $U(n/2)$, $SU(n/2)$, $Sp(n)$, $Sp(n)Sp(1)$,
$G_2$ and $Spin(7)$.
As is known to the physicists,
the $N=1$ SCVA in Theorem \ref{thm:N=1Riemannian} can be extended for
Riemannian manifolds of special holonomy groups. 
We will briefly indicate how this is done 
in the case of $G = G_2$ and $Spin(7)$ in this section
(cf. Shatashvilli-Vafa \cite{Sha-Vaf}).

Denote by $\rho: G_2 \to O(\bR^7)$ the holonomy representation of $G_2$,
and let $\{e^i: i = 1, \cdots, 7\}$ be an oriented orthonormal basis of $\bR^7$
with standard Euclidean metric $g_0$.
Then the $G_2$ preserves the following exterior form:
$$\Phi = e^1 \wedge e^2 \wedge e^5 + e^1 \wedge e^3 \wedge e^6 +
e^1 \wedge e^4 \wedge e^7 - e^2 \wedge e^3 \wedge e^7
+ e^2 \wedge e^4 \wedge e^6 - e^3 \wedge e^4 \wedge e^5
+ e^3 \wedge e^6 \wedge e^7.$$
Now regard $\Phi$ as an element of $\Lambda(t^{-\frac{1}{2}}\bC^7)$
in $V(\bC^7, g_0 \otimes \bC)$ and consider the field $\Phi(z)$.
More fields can be found by computing the OPE's.
For examples,
$$\Phi(z)\Phi(w) \sim \frac{-7}{(z-w)^3} + \frac{6X(w)}{z-w},$$
where 
\begin{eqnarray*}
X & = & - \phi^1_{-\frac{1}{2}} \wedge \phi_{-\frac{1}{2}}^2 
\wedge \phi_{-\frac{1}{2}}^3 \wedge \phi_{-\frac{1}{2}}^4 
+ \phi^1_{-\frac{1}{2}} \wedge \phi_{-\frac{1}{2}}^2 
\wedge \phi_{-\frac{1}{2}}^6 \wedge \phi_{-\frac{1}{2}}^7 \\
&& - \phi^1_{-\frac{1}{2}} \wedge \phi_{-\frac{1}{2}}^3 
\wedge \phi_{-\frac{1}{2}}^5 \wedge \phi_{-\frac{1}{2}}^7
+ \phi^1_{-\frac{1}{2}} \wedge \phi_{-\frac{1}{2}}^4 
\wedge \phi_{-\frac{1}{2}}^5 \wedge \phi_{-\frac{1}{2}}^6 \\
&& - \phi^2_{-\frac{1}{2}} \wedge \phi_{-\frac{1}{2}}^3 
\wedge \phi_{-\frac{1}{2}}^5 \wedge \phi_{-\frac{1}{2}}^6
- \phi^2_{-\frac{1}{2}} \wedge \phi_{-\frac{1}{2}}^4 
\wedge \phi_{-\frac{1}{2}}^5 \wedge \phi_{-\frac{1}{2}}^7 \\
&& - \phi^3_{-\frac{1}{2}} \wedge \phi_{-\frac{1}{2}}^4 
\wedge \phi_{-\frac{1}{2}}^6 \wedge \phi_{-\frac{1}{2}}^7
- \frac{1}{2} \phi_{-\frac{3}{2}}\phi_{-\frac{1}{2}}.
\end{eqnarray*}
Except for the last term, 
$X$ corresponds to $*\Phi$.
Two more states, $K$ and $M$, can be found by taking the OPE's
of $\Phi(z)$ and $X(z)$ with $G(z)$ respectively.
A remarkable fact is that the six states $\nu$, $\tau$, $\Phi$, $X$, $K$,
and $M$ generates a vertex algebra $V_{G_2}$ which is an extension of 
the $N=1$ Neveu-Schwarz algebra 
(cf. Shatashvilli-Vafa \cite{Sha-Vaf}, Appendix 1).
Furthermore,
it is easy to see that all these vectors are preserved by the $G_2$-action
induced by the holonomy representation.
Let $(M^7, g)$ be a Riemaniann manifold with holonomy group $G_2$.
Then there is a local trivialization of $TM$ such that the transition functions
lie in $G_2$. 
As a consequence,
each of these states induces a section in the bundle 
$V(TM \otimes \bC, g \otimes \bC)_{NS}$,
hence $V_{G_2}$ also acts on $\Gamma(M, V(TM \otimes \bC, g \otimes \bC)_{NS})$.
The case of $Spin(7)$ is similar.

\subsection{The case of quaternionic K\"{a}hler manifolds}

A {\em quaternionic K\"{a}hler manifold} of dimension $> 4$
is Riemannian manifold with holonomy group $Sp(n)Sp(1)$
(cf. Salamon \cite{Sal}).
Let $(M, g)$ be quaternionic K\"{a}hler manifold,
then there is a three dimensional subbundle $E \subset \Lambda^2(T^*M)$
preserved by the Levi-Civita connections,
such that locally one can find a basis $\{\omega_1, \omega_2, \omega_3\}$
corresponding to three local almost complex structures $K_1$, $K_2$, $K_3$ 
with $K_1K_2=-K_2K_1=K_3$.
The $4$-form
$$\Omega = \frac{1}{2} (\omega_1 \wedge \omega_1 
+ \omega_2 \wedge \omega_2 + \omega_3 \wedge \omega_3)$$
is globally well-defined and parallel.
Conversely a Riemannian manifold admiting a parallel $4$-form 
locally of this form is a quaternionic K\"{a}hler manifold.
Regard $\Omega$ as an element in 
$\Gamma(M, \Lambda^4(t^{-\frac{1}{2}}TM \otimes \bC)) \subset 
\Gamma(M, V(TM \otimes \bC, g \otimes \bC)_{NS})$.
It can also be regarded as an element $\widehat{\Omega}$ of 
$\Gamma(M, S^1(t^{-\frac{1}{2}}TM\otimes \bC) 
\otimes \Lambda^3(t^{-1}TM \otimes \bC))$.
Calculations by Wick's theorem in each fiber show that we have
\begin{eqnarray*}
L(z) \Omega(w) & \sim & \frac{2\Omega(w)}{(z-w)^2} 
	+ \frac{\partial_w\Omega(w)}{z-w}, \\
\Omega(z)\Omega(w) & \sim & \frac{3n(2n+1)}{(z-w)^4}
+ \frac{-4 \Omega(w) + 3n(2n+1)L_F(w)}{(z-w)^2} \\ 
&& + \frac{\partial_w(-4 \Omega(w) + 3n(2n+1) L_F(w))}{2(z-w)}, \\
G(z)\Omega(w) & \sim & \frac{\widehat{\Omega}(w)}{z-w}, \\
L(z)\widehat{\Omega}(w) & \sim & \frac{\frac{5}{2}\widehat{\Omega}(w)}{(z-w)^2}
+ \frac{\partial_w\widehat{\Omega}(w)}{z-w}, \\
G(z)\widehat{\Omega}(w) & \sim & \frac{4\Omega(w)}{(z-w)^2}
+\frac{\partial_w\Omega(w)}{z-w}, \\
\end{eqnarray*}
where $L_F(w)=Y(\nu_F,w)$.
The OPE's $\Omega(z)\widehat{\Omega}(w)$ and 
$\widehat{\Omega}(z)\widehat{\Omega}(w)$ are very complicated.
We will not pursue it here.

\section{$N=2$ SCVA bundles from complex manifolds}

\subsection{$N=2$ SCVA's from a space with an inner product and a polarization} 

A {\em polarization} of a vector space $T$ with an inner product
$g$ is a decomposition $T = T' \oplus T''$,
such that $g(a_1, a_2) = g(b_1, b_2) = 0$, 
for $a_1, a_2 \in T'$, $b_1, b_2 \in T''$.
Given a polarization $T = T' \oplus T''$ of $(T, g)$,
it is easy to see that $g$ induces an isomorphism $T'' \cong (T')^*$.
In particular,
this implies that a vector space with an inner product is even dimensional
if it admits a polarization.

From now on we assume that $(T, g)$ admits a polarization $T = T' \oplus T''$.
Then 
$$F(T, g)_{NS} = \otimes_{n>0} \Lambda(t^{-n + \frac{1}{2}}T') 
\otimes_{n > 0} \Lambda (t^{-n + \frac{1}{2}}T'').$$
Consider also the space
\begin{eqnarray*}
&& F(T, g)_R = \Lambda(\oplus_{n \geq 0} t^{-n} T') \otimes
\Lambda(\oplus_{n > 0} t^{-n} T'') \\
& = & \Lambda(T') \otimes_{n > 0} \Lambda(t^{-n}T') 
	\otimes_{n>0}\Lambda(t^{-n}T'')
\end{eqnarray*}
Notice that $\Lambda(T')$ is isomorphic to the space 
$\Delta(T)$ of spinors of $(T, g)$.
So we can also write
$$F(T, g)_R = \Delta(T) \otimes_{n > 0} \Lambda(t^{-n}T).$$
Let $\{\varphi^i\}$ be a basis of $T'$, 
$\{\psi^i\}$ a basis of $T''$,
such that $\omega(\varphi^i, \psi^j) = \delta_{ij}$.
Then $F(T, g)_R$ has a basis consists of elements of the form:
$$\varphi^{i_1}_{-k_1} \cdots \varphi^{i_m}_{-k_m}
\psi^{j_1}_{-l_1-1}\cdots\psi^{j_n}_{-l_n-1},$$
where $k_1, \cdots, k_m, l_1, \cdots, l_n \geq 0$.
In the NS case,
set
\begin{align*}
\varphi^i(z) & = \sum_{r \in \frac{1}{2} + \bZ}
	\varphi^i_{(r)} z^{-r - \frac{1}{2}}, &
\psi^i(z) & = \sum_{r \in \frac{1}{2} +\bZ} \psi^i_{(r)} z^{-r - \frac{1}{2}},
\end{align*}
where for $r < 0$, 
$\varphi^i_{(r)}$ and $\psi^i_{(r)}$ are exterior products by
$\varphi^i_r$ and $\psi^i_r$ respectively,
for $ r > 0$,
$\varphi^i_{(r)}$ and $\psi^i_{(r)}$ are contractions by
$\varphi^i_{-r}$ and $\psi^i_{-r}$ respectively;
in the R case,
set
\begin{align*}
\varphi^i(z) & = \sum_{n \in \bZ} \varphi^i_n z^{-n}, &
\psi^i(z) & = \sum_{n \in \bZ} \psi^i_n z^{-n-1},
\end{align*}
where for $n < 0$, 

Similar to Proposition \ref{prop:N=1},
we have the following the following:

\begin{proposition} \label{prop:N=2}
Let $(T, g)$ be a finite dimensional complex vector space
with a polarization $T = T' \oplus T''$.
Also let $\{e^i\}$ and $\{f^j\}$ be bases of $T'$ and $T''$
such that $g(e^i, f^j) = \delta_{ij}$.
Write these bases as $\{b^i\}$ and $\{c^j\}$ respective 
for the copies of $T$ in the bosonic sector,
and as $\{\varphi^i\}$ and $\{\psi^i\}$ for the copies of $T$ 
in the fermionic sector.

(a) There is a structure of a vertex algebra on
$F(T, g)_{NS}$ defined by
\begin{eqnarray*}
&& Y(\varphi^{i_1}_{-k_1 - \frac{1}{2}} \cdots \varphi^{i_m}_{-k_m-\frac{1}{2}}
\psi^{j_1}_{-l_1 -\frac{1}{2}} \cdots \psi^{j_n}_{-l_n-\frac{1}{2}}) \\
& = & 
:\partial_z^{(k_1)}\varphi^{i_1}(z) \cdots \partial_z^{(k_m)}\varphi^{i_m}(z)
\partial_z^{(l_1)}\psi^{j_1}(z)\cdots\partial_z^{(l_m)}\psi^{j_m}(z):.
\end{eqnarray*}
Furthermore,
for any $\lambda \in \bC$,
$$\nu_{\lambda} 
= (1-\lambda)\sum_i \varphi^i_{-\frac{3}{2}}\psi^i_{-\frac{1}{2}}
+ \lambda \sum_i \psi^i_{-\frac{3}{2}}\varphi^i_{-\frac{1}{2}}$$ 
is a conformal vector of central charge 
$$-(6 \lambda^2 - 6 \lambda + 1) \dim T.$$
Simlarly for $F(T, g)_R$.

(b) There is a structure of is an $N=2$ SCVA of conformal weight 
$\frac{3}{2} \dim T$
on  $V(T, g)_{NS} = F(T, g)_{NS} \otimes B(T, g)$ 
with superconformal structures given by the following vectors:
\begin{align*}
\tau^+ & = \sum_i b_{-1}^i\psi^i_{-\frac{1}{2}}, &
\tau^- & = \sum_i c^i_{-1}\varphi^i_{-\frac{1}{2}}, \\
j & = \sum_i \psi^i_{-\frac{1}{2}}\varphi^i_{-\frac{1}{2}}, &
\nu & = \sum_i (b^i_{-1}c^i_{-1} 
	+ \frac{1}{2}\varphi^i_{-\frac{3}{2}}\psi^i_{-\frac{1}{2}}
	+ \frac{1}{2}\psi^i_{-\frac{3}{2}}\varphi^i_{-\frac{1}{2}}).
\end{align*}
Similarly for $V(T, g)_R = F(T, g)_R \otimes B(T, g)$
with 
$\varphi^j_{-\frac{1}{2}}$ and $\varphi^j_{-\frac{3}{2}}$ replaced by
$\varphi^j_0$ and $\varphi^j_{-1}$ respectively,
and $\psi^j_{-\frac{1}{2}}$ and $\psi^j_{-\frac{3}{2}}$ replaced by
$\psi^j_{-1}$ and $\psi^j_{-2}$ respectively.
\end{proposition}

\begin{remark} \label{rmk:N=2}
It is clear that all the vectors in Proposition \ref{prop:N=2}
are independent of the choice of the bases. 
Denote by $GL(T')$ the group of linear transformation on $T'$.
Since $T'' \cong (T')^*$,
there is an induced action of $GL(T')$ on $T''$,
hence $GL(T')$ acts on $T$, preserving $g$.
This action extends to actions on $F(T,g)_{NS}$ and $F(T, g)_R$
as automorphisms of vertex algebras,
and to actions on $V(T, g)_{NS}$ and $V(T, g)_R$ as automorphisms
of $N=2$ SCVA's.
\end{remark}

For later use,
we remark that in $V(T' \oplus T'', g)_R$, 
\begin{align*}
& L_0b^i_{-n} = n b^i_{-1}, & L_0 c^i_{-1} & = n c^i_{-n}, 
& L_0\varphi^i_{-n + 1} = (n - \frac{1}{2}) \varphi^i_0, & 
& L_0\psi^i_{-n} = (n - \frac{1}{2}) \psi^i_{-1}, \\
& J_0 b^i_{-1} = 0, & J_0 c^i_{-1} & = 0, 
& J_0\varphi^i_{- n + 1} = -\varphi^i_0, & 
& J_0\psi^i_{-n} = \psi^i_{-1}, 
\end{align*} 
hence
\begin{eqnarray} \label{elliptic:N=2}
&& G_{q, y}(V(T, g)_R) \\ 
& = & q^{-\frac{\dim T}{16}}
\otimes_{n > 0} \Lambda_{y^{-1}q^{n-\frac{1}{2}}}(T')
\otimes_{n > 0} \Lambda_{yq^{n-\frac{1}{2}}}(T'')
\otimes_{n > 0} S_{q^n}(T') \otimes_{n > 0} S_{q^n}(T''). \notag
\end{eqnarray}
In the $A$ twist we have
\begin{align*}
& T_0b^i_{-n} = n b^i_{-1}, & T_0 c^i_{-1} & = n c^i_{-n}, 
& T_0\varphi^i_{-n + 1} = n \varphi^i_0, & 
& T_0\psi^i_{-n} = (n - 1) \psi^i_{-1}, \\
& J^{top}_0 b^i_{-1} = 0, & J^{top}_0 c^i_{-1} & = 0, 
& J^{top}_0\varphi^i_{- n + 1} = -\varphi^i_0, & 
& J^{top}_0\psi^i_{-n} = \psi^i_{-1}, 
\end{align*} 
hence
\begin{eqnarray}  \label{elliptic:Atwist}
&& G_{q, y}(V(T, g)_R) \\
& = & q^{-\frac{\dim T}{16}}
\otimes_{n > 0} \Lambda_{y^{-1}q^{n}}(T')
\otimes_{n > 0} \Lambda_{yq^{n-1}}(T'')
\otimes_{n > 0} S_{q^n}(T') \otimes_{n > 0} S_{q^n}(T''). \notag
\end{eqnarray}
Similarly,
in the $B$ twist we have
\begin{align*}
& T_0b^i_{-n} = n b^i_{-1}, & T_0 c^i_{-1} & = n c^i_{-n}, 
& T_0\varphi^i_{-n + 1} = (n - 1) \varphi^i_0, & 
& T_0\psi^i_{-n} = n \psi^i_{-1}, \\
& J^{top}_0 b^i_{-1} = 0, & J^{top}_0 c^i_{-1} & = 0, 
& J^{top}_0\varphi^i_{- n + 1} = \varphi^i_0, & 
& J^{top}_0\psi^i_{-n} = - \psi^i_{-1}, 
\end{align*} 
hence
\begin{eqnarray}\label{elliptic:Btwist}
&& G_{q, y}(V(T, g)_R) \\
& = & q^{-\frac{\dim T}{16}}
\otimes_{n > 0} \Lambda_{yq^{n-1}}(T')
\otimes_{n > 0} \Lambda_{y^{-1}q^{n}}(T'')
\otimes_{n > 0} S_{q^n}(T') \otimes_{n > 0} S_{q^n}(T''). \notag
\end{eqnarray}

Similar to Theorem 2.4 in Malikov-Schechtman-Vaintrob \cite{Mal-Sch-Vai},
we have the following:

\begin{theorem} \label{thm:BRSTinner}
(a) For the $N=2$ SCVA $V(T' \oplus T'', g)_{NS}$ 
in Proposition \ref{prop:N=2},
the topological vertex algebras obtained by $A$ twist and $B$ twist 
(cf. Proposition \ref{prop:top}) both have trivial BRST cohomology.

(b) For the $N=2$ SCVA $V(T' \oplus T'', g)_R$ in Proposition \ref{prop:N=2},
the BRST cohomology of the topological vertex algebras obtained by
$A$ twist and $B$ twist (cf. Proposition \ref{prop:top})
are isomorphic to $\Lambda(T'')$ and $\Lambda(T')$ 
as graded commutative algebras respectively.
\end{theorem}

\begin{proof}
(a) Let $Q_0$ be the zero mode of $Q(z) = G^+(z)$. 
We have
$$Q = \sum_{n < 0} \sum_i b_n^i\psi^i_{-n-\frac{1}{2}}
+ \sum_{n \geq 0} \sum_i \psi^i_{-n-\frac{1}{2}}b_n^i.$$
Since $b_0^i$ acts as $0$,
we actually have $Q_0 = Q_- + Q_+$,
where
\begin{align*}
Q_- & = \sum_{n < 0} \sum_i b_n^i\psi^i_{-n-\frac{1}{2}}, &
Q_+ & = \sum_{n > 0} \sum_i \psi^i_{-n-\frac{1}{2}}b_n^i.
\end{align*}
It is easy to see that $Q_-^2 = [Q_-, Q_+] = Q_+^2 = 0$,
hence we get a double complex and two spectral sequences
with $E_1$ term the $Q_+$-cohomology and the $Q_-$-cohomology 
respectively (cf. Bott-Tu \cite{Bot-Tu}). 
Now 
\begin{eqnarray*}
&& V_{NS}(T, g) \\
& = & \bigotimes_{n>0} \left(\Lambda(t^{-n+\frac{1}{2}}T') 
\otimes S(t^{-n}T'') \right) \otimes \bigotimes_{n>0} \left(
\Lambda(t^{-n+\frac{1}{2}}T'') \otimes S(t^{-n}T') \right).
\end{eqnarray*}
On the first factor,
$Q_-$ acts as the differential in 
the tensor product of infinitely many copies of Koszul complexes,
while $Q_+$ acts trivially;
on the second factor,
$Q_-$ acts trivially,
while $Q_+$ acts as the differential in the tensor product of 
infinitely many copies of algebraic de Rham complexes of $V$.
The proof is completed by taking cohomology in $Q_-$ then in $Q_+$.
The case of $Q(z) = G^-(z)$ is similar.

(b) Similarly in the case of $Q(z) = G^+(z)$,
we have
$Q_0 = Q_- + Q_+$,
where
\begin{align*}
Q_- & = \sum_{n < 0} \sum_i b_n^i\psi^i_{-n-1}, &
Q_+ & = \sum_{n > 0} \sum_i \psi^i_{-n-1}b_n^i.
\end{align*}
Again we have $Q_-^2 = [Q_-, Q_+] = Q_+^2 = 0$. 
Now 
\begin{eqnarray*}
V_R(T, g) & = & \bigotimes_{n>0} \left(\Lambda(t^{-n+1}T') 
\otimes S(t^{-n}T'') \right) \\
&& \otimes \bigotimes_{n>0} \left(
\Lambda(t^{-n-1}T'') \otimes S(t^{-n}T') \right) \otimes \Lambda(t^{-1}T'').
\end{eqnarray*}
Similar to (a),
one sees that cohomology in $Q$ is $\Lambda(t^{-1}T'')$.
Now any element of $\Lambda(t^{-1}T'')$ is of the form
$$\psi^{j_1}_{-1} \cdots \psi^{j_n}_{-1}$$
for some $j_1, \cdots, j_n$.
It corresponds to the field
$$:\psi^{j_1}(z) \cdots \psi^{j_n}(z):.$$
Given two elements $\psi^{j_1}_{-1} \cdots \psi^{j_n}_{-1}$
and $\psi^{k_1}_{-1} \cdots \psi^{k_m}_{-1}$,
by Wick's theorem,
we have
\begin{eqnarray*}
&& :(:\psi^{j_1}(z) \cdots \psi^{j_n}(z):)
(:\psi^{k_1}(z) \cdots \psi^{k_n}(z):): \\
& = & :\psi^{j_1}(z) \cdots \psi^{j_n}(z)\psi^{k_1}(z) \cdots \psi^{k_n}(z):.
\end{eqnarray*}
Hence on the $Q_0$-cohomology, 
the product induced from the normally ordered product is isomorphic to
the ordinary exterior product on $\Lambda(t^{-1}T'')$.
The case of $Q(z) = G^-(Z)$ is similar.
\end{proof}

\subsection{$N=2$ SCVA bundles from complex manifolds}
\label{sec:complex}

Let $(M, J)$ be a complex manifold.
Denote by $T_cM$ the holomorphic tangent bundle.
The fiberwise pairing between $T_cM$ and $T_c^*M$
induces a canonical complex inner product $\eta$ on 
the holomorphic vector bundle $T_cM \oplus T_c^*M$
with a manifest polarization $T' = T_cM$, $T''=T_c^*M$.
By construction of Proposition \ref{prop:N=2} and Remark \ref{rmk:N=2},
we obtain an $N=2$ SCVA bundle $V(T_cM \oplus T_c^*M, \eta)_R$.
Since this bundle is holomorphic,
one can consider the $\bar{\partial}$ operator on it:
$$\bar{\partial}: \Omega^{0, *}(V(T_cM \oplus T_c^*M, \eta)_R) 
\to \Omega^{0, *+1}(V(T_cM \oplus T_c^*M, \eta)_R).$$

\begin{theorem} \label{thm:N=2complex}
For any complex manifold $M$,
$\Omega^{0, *}(V(T_cM \oplus T_c^*M, \eta)_R)$ has a natural structure of
an $N=2$ SCVA such that $\bar{\partial}$ is a differential.
Consequently,
the Dolbeault cohomology 
$$H^*(M, V(T_cM \oplus T_c^*M, \eta)_R)$$
is an $N=2$ SCVA;
furthermore,
the BRST cohomology of
its associated topological vertex algebras (cf. Proposition \ref{prop:top})
is isomorphic
to $H^*(M, \Lambda(T_cM)$ or $H^*(M, \Lambda(T^*_cM))$ depending on
whether we take $Q(z) = G^+(z)$ or $G^-(z)$.
Similar results can be obtained for $V(T_cM \oplus T_c^*M, \eta)_{NS}$.
However the BRST cohomologies are trivial 
for the corresponding Dolbeault cohomology.
\end{theorem}

\begin{proof}
We regard $\Lambda(\overline{T_c}^*M)$ as a bundle of holomorphic 
vertex algebra,
therefore, 
$V(T_cM \oplus T_c^*M, \eta)_R \otimes \Lambda(\overline{T_c}^*M)$
has a natural structure of an $N=2$ SCVA.
By Lemma \ref{lm:section},
the section space
$$\Gamma(M, V(T_cM \oplus T_c^*M, \eta)_R 
\otimes \Lambda(\overline{T_c}^*M))$$
is an $N=2$ SCVA.
One can easily verify that $\bar{\partial}$ is a differential by
choosing a local holomorphic frame of $T_cM$.
It follows from Lemma \ref{lm:cohomology} that
$H^*(M, V(T_cM \oplus T_c^*M, \eta)_R)$
is an $N=2$ SCVA.
Notice that on $\Gamma(M, V(T_cM \oplus T_c^*M, \eta)_R 
\otimes \Lambda(\overline{T_c}^*M))$
two operators $\bar{\partial}$ and $Q_0$ act such that
$$\bar{\partial}^2 = [\bar{\partial}, Q_0] = Q_0^2 = 0.$$
In the above we have taken the $\bar{\partial}$-cohomology first,
then take the $Q_0$-cohomology.
We can also do it in a different order.
By Theorem \ref{thm:BRSTinner},
the $Q_0$-cohomology is $\Lambda(T_cM)$ or $\Lambda(T^*_cM)$,
its $\bar{\partial}$-cohomology is the Dolbeault cohomology.
This completes the proof.
\end{proof}

As in $N=1$ vertex algebra bundle in the Riemannian case 
(cf. \S \ref{sec:Riemannian}),
$V(T_cM \oplus T_c^*M, \eta)_R$ is related to the elliptic genera.
By (\ref{elliptic:N=2}) we have
\begin{eqnarray*}
&& G_{q, y}(V(T_cM \oplus T_c^*M, \eta)_R) \\
& = & \otimes_{n > 0} \Lambda_{y^{-1}q^{n-\frac{1}{2}}}(T_cM) 
\otimes_{n > 0} \Lambda_{yq^{n-\frac{1}{2}}}(T_c^*M) 
\otimes_{n > 0} S_{q^n}(T_cM) \otimes_{n > 0} S_{q^n}(T_c^*M).
\end{eqnarray*}
By (\ref{elliptic:Atwist}) we see that in the $A$ twist we have
\begin{eqnarray*}
&& G_{q, y}(V(T_cM \oplus T_c^*M, \eta)_R) \\
& = & \otimes_{n > 0} \Lambda_{y^{-1}q^n}(T_cM) 
\otimes_{n > 0} \Lambda_{yq^{n-1}}(T_c^*M) 
\otimes_{n > 0} S_{q^n}(T_cM) \otimes_{n > 0} S_{q^n}(T_c^*M)
\end{eqnarray*}
(cf. Hirzebruch \cite{Hir}, (16)),
while in the $B$ twist we have
\begin{eqnarray*}
&& G_{q, y}(V(T_cM \oplus T_c^*M, \eta)_R) \\
& = & \otimes_{n > 0} \Lambda_{yq^{n-1}}(T_cM) 
\otimes_{n > 0} \Lambda_{y^{-1}q^n}(T_c^*M) 
\otimes_{n > 0} S_{q^n}(T_cM) \otimes_{n > 0} S_{q^n}(T_c^*M)
\end{eqnarray*}
(cf. Dijkgraaf {\em et.~al.} \cite{Dij-Moo-Ver-Ver}, (A.8)).

\subsection{$N=2$ SCVA bundles from K\"{a}hler manifolds}
\label{sec:Kahler}

Assume that $(M, J, g, \omega)$ is a K\"{a}hler manifold.
In the decomposition $TM \otimes \bC = T'M \oplus T''M$,
where $T'M \cong T_cM$,
$T''M \cong \overline{T'}M \cong \overline{T_c}M$.
The K\"{a}hler metric induces an isomorphism of complex vector bundles:
$T''M \cong T_c^*M$.
Under this isomorphism,
the complexified metric $g_{\bC}$ can be naturally identified with the
complex canonical inner product $\eta$ on $T_cM \oplus T^*_cM$ used above.
Hence we have
$$V(TM \otimes \bC, g\otimes \bC)_R 
\cong V(T_cM \oplus T^*_cM, \eta)_R.$$
Using the $(0,1)$-part of the Levi-Civita connection,
a $\bar{\partial}$ operator can be defined on
$$\Omega^{0, *}(V(TM \otimes \bC, g\otimes \bC)).$$
It is naturally identified with the $\bar{\partial}$ operator on 
$\Omega^{0, *}(V(T_cM \oplus T_c^*M, \eta)_R)$.
Hence we have the following analogue of Theorem \ref{thm:N=2complex} for 
$V(TM \otimes \bC, g\otimes \bC)_R$ in the case of K\"{a}hler manifolds:

\begin{theorem} \label{thm:N=2Kahler}
For any K\"{a}hler manifold $M$,
$\Omega^{0, *}(V(TM \otimes \bC, g \otimes \bC)_R)$ 
has a natural structure of an $N=2$ SCVA and a differential $\bar{\partial}$.
Consequently,
the  cohomology 
$$H^*(\Gamma(M, V(TM \otimes \bC, g \otimes \bC)_R), \bar{\partial})$$
is an $N=2$ SCVA;
furthermore,
the BRST cohomology of
its associated topological vertex algebras (cf. Proposition \ref{prop:top})
is isomorphic
to $H^*(M, \Lambda(T_cM)$ or $H^*(M, \Lambda(T^*_cM))$ depending on
whether we take $Q(z) = G^+(z)$ or $G^-(z)$.
Similar results can be obtained for $V(TM \otimes \bC, g \otimes \bC)_{NS}$.
However the BRST cohomologies are trivial 
for the corresponding Dolbeault cohomology.
\end{theorem}

\subsection{Extended $N=2$ SCVA's from Calabi-Yau manifolds}
\label{sec:CY}

When the K\"{a}hler manifold $(M^n, g, J)$ is Calabi-Yau 
and has holonomy group $SU(n)$,
then the $N=2$ supersymmetry of the $N=2$ SCVA's 
constructed in \S \ref{sec:Kahler} can be extended 
to algebras described in Odake \cite{Oda} 
(see also Figueroa-O'Farill \cite{Fig}).
When $n=2$, 
this manifold is also a hyperk\"{a}hler manifold 
for which we will discuss the $N=4$ supersymmetry below.
Recall that there is a nonzero parallel holomorphic $n$-form $\Omega$ on $M$.
Without loss of generality,
we can take local parallel orthonormal frame $\{e^i\}$ of $(T'M)^*$,
such that
\begin{align*}
\omega & = \sum_i e^i \wedge \bar{e}^i, &
\Omega & = e^1 \wedge \cdots\wedge e^n.
\end{align*}
Let
\begin{align*}
X^+(z) & = :\psi^1(z)\cdots\psi^n(z):, &
Y^+(z) & = \sum_{j=1}^n (-1)^{j-1} :c^j(z)\psi^1(z)\cdots 
\widehat{\psi^j(z)}\cdots \psi^n(z):, \\
X^-(z) & = :\varphi^1(z)\cdots\varphi^n(z):, &
Y^-(z) & = \sum_{j=1}^n (-1)^{j-1} :b^j(z)\varphi^1(z)\cdots 
\widehat{\varphi^j(z)}\cdots \varphi^n(z):.
\end{align*}
Then $X^{\pm}(z)$ and $Y^{\pm}(z)$ are globally well-defined sections of
$V(TM \otimes \bC, g\otimes\bC)_{NS}$ or $V(TM \otimes \bC, g\otimes \bC)_R$.
By fiberwise calculations by Wick's theorem,
we get
\begin{eqnarray*}
&& L(z)X^{\pm}(w) \sim \frac{nX^{\pm}(w)}{2(z-w)^2}, \;\;\;\;
J(z)X^{\pm}(w) \sim \pm \frac{nX(w)}{z-w}, \\
&& L(z)Y^{\pm}(w) \sim \frac{(n+1)Y^{\pm}(w)}{2(z-w)^2}, \;\;\;\;
J(z)Y^{\pm}(w) \sim \pm \frac{(n-1)Y(w)}{z-w}, \\
&& G^+(z)X^+(w) \sim 0, \;\;\;\; G^-(z)X^+(w) \sim \frac{Y^+(w)}{z-w}, \\
&& G^+(z)Y^+(w) \sim \frac{nX^+(w)}{(z-w)^2} + \frac{\partial_wX^+(w)}{z-w}, 
\;\;\;\; G^-(z)Y^+(w) \sim 0, \\
&& G^+(z)X^-(w) \sim Y^-(w), \;\;\;\; G^-(z)X^+(w) \sim 0, \\
&& G^+(z)Y^-(w) \sim 0,  \;\;\;\; 
G^-(z)Y^+(w) \sim \frac{nX^+(w)}{(z-w)^2} + \frac{\partial_wX^-(w)}{z-w}, \\
&& X^{\pm}(z)X^{\pm}(w) \sim X^{\pm}(z)Y^{\pm}(w) 
\sim Y^{\pm}(z)Y^{\pm}(w) \sim 0.
\end{eqnarray*}
Hence $\{L(z), J(z), G^{\pm}(z), X^+(z), Y^+(z)\}$ generate a vertex subalgebra.
Similarly for $\{L(z), J(z), G^{\pm}(z), X^-(z), Y^-(z)\}$.
When 
When $n=2$,
we also have
\begin{eqnarray*}
&& X^+(z)X^-(w) \sim -\frac{1}{(z-w)^2} - \frac{J(w)}{z-w},\\
&& X^+(z)Y^-(w) \sim \frac{G^+(w)}{z-w}, \;\;\;\;
X^-(z)Y^+(w) \sim  \frac{G^-(w)}{z-w}, \\
&& Y^+(z)Y^-(w) \sim \frac{2}{(z-w)^3} + \frac{J(w)}{(z-w)^2} 
+ \frac{L(w) + \frac{1}{2}\partial_wJ(w)}{z-w}.
\end{eqnarray*}
Hence $\{L(z), J(z), G^{\pm}(z), X^{\pm}(z), Y^{\pm}(z)\}$
generate a vertex subalgebra.
When $n=3$,
we have
\begin{eqnarray*}
X^+(z)X^-(w) & \sim & -\frac{1}{(z-w)^3} - \frac{J(w)}{(z-w)^2}
- \frac{:J(w)J(w): - \partial_wJ(w)}{2(z-w)},\\
X^+(z)Y^-(w)  & \sim & -\frac{G^+(w)}{(z-w)^2} - \frac{:J(w)G^+(w):}{z-w}, \\
X^-(z)Y^+(w)  & \sim & -\frac{G^-(w)}{(z-w)^2} + \frac{:J(w)G^-(w):}{z-w}, \\
Y^+(z)Y^-(w)  & \sim & \frac{-3}{(z-w)^4} 
- \frac{2J(w)}{(z-w)^3}
- \frac{\frac{1}{2}:J(w)J(w): + L(w) - \partial_w J(w)}{(z-w)^2}.
\end{eqnarray*}
The calculations for bigger $n$ is similar but more complicated. 
(See \cite{Fig} for results when $n=4$.)

\section{$N=4$ SCVA's from hyperk\"{a}hler manifolds}
\label{sec:N=4}

\subsection{Definition of an $N=4$ SCVA}

The ``small" $N=4$ superconformal algebra was introduced by 
Ademollo {\em et.~al.} \cite{Ade}.
Here we adapt the notations in Berkovits and Vafa \cite{Ber-Vaf}.
Recall that the small $N=4$ superconformal algebra 
is the $N=2$ superconformal algebra $(L, G^{\pm}, J)$ with two additional
currents $J^{++}$ and $J^{--}$ of charge $\pm 2$,
and two supercurrents $\tilde{G}^+$ and $\tilde{G}^-$.
Besides the OPE's of the $N=2$ superconformal algebra,
the following additional OPE's are required:
\begin{align*}
& L(z)\tilde{G}^{\pm}(w) \sim \frac{3}{2} \frac{\tilde{G}^{\pm}(w)}{(z-w)^2}
 + \frac{\partial_w\tilde{G}^{\pm}(w)}{z-w}, \;\;\;\;
J(z) \tilde{G}^{\pm}(w) \sim \frac{\pm \tilde{G}(w)}{z-w}, \\
& L(z) J^{\pm\pm}(w) \sim \frac{J^{\pm\pm}(w)}{(z-w)^2} 
+ \frac{\partial_wJ^{\pm}(w)}{z-w}, \;\;\;\;
J(z)J^{\pm}(w) \sim \frac{\pm 2 J^{\pm}(w)}{z-w},  \\
& J^{--}(z)J^{++}(w) \sim - \frac{c}{3(z-w)^2} + \frac{J(w)}{z-w}, \;\;\;\;
J^{\pm}(z) \cdot J^{\pm}(w) \sim 0, \\
& J^{--}(z)G^{+}(w) \sim \frac{\tilde{G}^-(w)}{z-w}, \;\;\;\;
J^{++}(z)\tilde{G}^{-}(w) \sim \frac{-G^+(w)}{z-w}, \\
& J^{++}(z)G^{-}(w) \sim \frac{\tilde{G}^+(w)}{z-w}, \;\;\;\;
J^{--}(z)\tilde{G}^{+}(w) \sim \frac{-G^-(w)}{z-w}, \\
& J^{--}(z)G^-(w) \sim J^{++}(z)G^+(w) \sim J^{++}(z)\tilde{G}^+(w) \sim 
J^{--}(z)\tilde{G}^-(w) \sim 0, \\
& G^+(z)\tilde{G}^-(w) \sim G^-(z)\tilde{G}^+(w) \sim 0, \\
& \tilde{G}^+(z)\tilde{G}^-(w) \sim \frac{2c}{3(z-w)^3}
	+ \frac{J(w)}{(z-w)^2} 
	+ \frac{L(w) + \frac{1}{2} \partial_wJ(w)}{z-w}, \\
& G^+(z)\tilde{G}^+(w) \sim -\frac{2J^{++}(w)}{(z-w)^2} 
	- \frac{\partial_wJ^{++}(w)}{z-w}, \\
& G^-(z)\tilde{G}^-(w) \sim -\frac{2J^{--}(w)}{(z-w)^2} 
	- \frac{\partial_wJ^{--}(w)}{z-w}, \\
\end{align*}
An {\em $N=4$ superconformal vertex algebra} is a 
vertex algebra $V$ with vectors $\nu$, $\tau^{\pm}$, $\tilde{\tau}^{\pm}$,
$j$, $j^{++}$ and $j^{--}$
such that $L(z) = Y(\nu, z)$, $G^{\pm}(z) = Y(\tau^{\pm}, z)$, 
$\tilde{G}^{\pm}(z) = Y(\tilde{\tau}^{\pm}, z)$,
$J(z) = Y(j, z)$, $J^{++}(z) = Y(j^{++}, z)$, 
and $J^{--}(z) = Y(j^{--}, z)$ generates an $N=4$ superconformal algebra.

\subsection{$N=4$ SCVA from a quternionic vector space}
Let $(T, g)$ be a $4n$-dimensional real vector space 
with three almost complex structures 
$K_1$, $K_2$ and $K_3$ which are compatible with metric and such that 
$K_1K_2 = - K_2K_1 = K_3$.
Such a space will be referred to as a {\em quaternionic vector space}.
Consider the polarization $T \otimes \bC = T' \oplus T''$,
where $K_1|_{T'} = \sqrt{-1}$, $K_1|_{T''} = - \sqrt{-1}$.
Choose an orthonormal basis of $T$ of the form
$\{a^i, b^i, c^i, d^i\}$,
where $b^i = K_1a^i$, $c^i = K_2a^i$, $d^i = K_3 a^i$.
Let 
\begin{align*}
e^{2i-1} & = \frac{1}{\sqrt{2}}(a^i - \sqrt{-1} b^i), & 
e^{2i} & = \frac{1}{\sqrt{2}}(c^i - \sqrt{-1} d^i), \\
\bar{e}^{2i-1} & = \frac{1}{\sqrt{2}}(a^i + \sqrt{-1} b^i), & 
\bar{e}^{2i} & = \frac{1}{\sqrt{2}}(c^i + \sqrt{-1} d^i).
\end{align*}
Then $\{e^i\}$ is a basis of $T'$, 
and $\{\bar{e}^i\}$ is a basis of $T''$,
such that $g(e^i, e^j) = g(\bar{e}^i, \bar{e}^j) = 0$,
$g(e^i, \bar{e}^j) = \delta_{ij}$.
As usual, 
these bases will be written as $\{b^i\}$ and $\{c^i\}$ for the copies
in the bosonic sector,
as $\{\varphi^i\}$ and $\{\psi^i\}$ for the copies in the fermionic sector.

\begin{proposition} \label{prop:N=4}
Let $(T, g, K_1, K_2, K_3)$ be as above.
Then $V(T,g)_{NS}$ has a structure of an $N=4$ SCVA of central charge 
$\frac{3}{2}\dim T$ given by the following vectors:
\begin{eqnarray*}
&& \nu = \sum_{i=1}^{\dim T/2}(b^i_{-1}c^i_{-1} 
+ \frac{1}{2} \varphi^i_{-\frac{3}{2}} \psi^i_{-\frac{1}{2}}
+ \frac{1}{2} \psi^i_{-\frac{3}{2}} \varphi^i_{-\frac{1}{2}}), \\
&& \tau^+ = \sum_{i=1}^{\dim T/2} b^i_{-1}\psi^i_{-\frac{1}{2}}, \;\;\;\;
\tau^- = \sum_{i=1}^{\dim T/4} c^i_{-1}\varphi^i_{-\frac{1}{2}}, \\
&& \tilde{\tau}^+ = \sum_{i=1}^{\dim T/4} 
	(c^{2i-1}_{-1}\psi^{2i}_{-\frac{1}{2}} 
	- c^{2i}_{-1}\psi^{2i-1}_{-\frac{1}{2}}), \;\;\;\;
\tilde{\tau}^-\sum_{i=1}^{\dim T/4} 
	(b^{2i-1}_{-1}\varphi^{2i}_{-\frac{1}{2}} 
	- b^{2i}_{-1}\varphi^{2i-1}_{-\frac{1}{2}}), \\
&& j = \sum_{i=1}^{\dim T/2} \psi^i_{-\frac{1}{2}}\phi^i_{-\frac{1}{2}}, 
\;\;\;\;
j^{++} = \sum_{i=1}^{\dim T/4} 
\psi^{2i}_{-\frac{1}{2}}\psi^{2i-1}_{-\frac{1}{2}}, \;\;\;\;
j^{--} = \sum_{i=1}^{\dim T/4} 
\varphi^{2i}_{-\frac{1}{2}}\varphi^{2i-1}_{-\frac{1}{2}}.
\end{eqnarray*}  
Similar result holds for $V(T, g)_R$ with 
$\varphi^j_{-\frac{1}{2}}$ and $\varphi^j_{-\frac{3}{2}}$ replaced by
$\varphi^j_0$ and $\varphi^j_{-1}$ respectively,
and $\psi^j_{-\frac{1}{2}}$ and $\psi^j_{-\frac{3}{2}}$ replaced by
$\psi^j_{-1}$ and $\psi^j_{-2}$ respectively.
\end{proposition}

\begin{remark} \label{rmk:N=4}
It is easy to see that the vectors in Proposition \ref{prop:N=4}
are independent of the choice of the basis.
Furthermore,
$Sp(n)$ acts by auotmorphisms of $N=4$ SCVA.
\end{remark}

\subsection{$N=4$ SCVA bundles from hyperk\"{a}hler manifolds}
\label{sec:hyperkahler}

We now generalize the results in \S \ref{sec:Kahler}.
Let $(M^{4n}, g, K_1, K_2, K_3)$ be a hyperK\"{a}hler manifold.
Then it is possible to find trivializations of $TM$ 
such that all the transition functions lie in $Sp(n)$.
By construction of Proposition \ref{prop:N=4} and Remark \ref{rmk:N=4},
we obtain an $N=4$ SCVA bundle $V(TM \otimes \bC, g\otimes \bC)_R$.
Similar to the argument in \S \ref{sec:Kahler},
one can consider the $\bar{\partial}_{K_1}$ operator on it:
$$\bar{\partial}_{K_1}: \Omega^{0, *}(V(TM \otimes \bC, g\otimes \bC)_R) 
\to \Omega^{0, *+1}(V(TM \otimes \bC, g\otimes \bC)_R).$$

\begin{theorem} \label{thm:N=4hyperkahler}
For any hyperk\"{a}hler manifold $M$,
$\Omega^{0, *}(V(TM \otimes \bC, g\otimes \bC)_R)$ has a natural structure of
an $N=4$ SCVA such that $\bar{\partial}_{K_1}$ is a differential.
Consequently,
the Dolbeault cohomology 
$$H^*(M, V(TM \otimes \bC, g\otimes \bC)_R)$$
is an $N=4$ SCVA;
furthermore,
the BRST cohomology of
its associated topological vertex algebras (cf. Proposition \ref{prop:top})
is isomorphic
to $H^*(M, \Lambda(T_cM))$ or $H^*(M, \Lambda(T^*_cM))$ depending on
whether we take $Q(z) = G^+(z)$ or $G^-(z)$.
Similar results can be obtained for $V(TM \otimes \bC, g\otimes \bC)_{NS}$.
However the BRST cohomologies are trivial for 
the corresponding Dolbeault cohomology.
\end{theorem}

\section{Relationship with loop spaces}

Our work is inspired by Vafa's suggestion \cite{Vaf} 
of an approach to quantum cohomology
based on vertex algebra constructed via semi-infinite forms on loop space.
Recall that a closed string in a manifold $M$ is a smooth map 
from $S^1$ to $M$.
The space of all closed string is just the free loop space $LM$.
Though extensively studied in algebraic topology,
earlier researches mostly dealt with the ordinary cohomology of the loop spaces.
Experience with infinite dimensional algebra and geometry 
(cf. Feigen-Frenkel \cite{Fei-Fre} and Atiyah \cite{Ati2})
tells us that we should look at the cohomology theory related to
semi-infinite forms on the loop space instead.
The space of such forms has a natural structure of a vertex algebra 
being the Fock space of a natural infinite dimensional Clifford algebra.
Vafa suggested to look at the multiplication in this vertex algebra.
We regard this multiplication as an infinite generalization of
the Clifford multiplication and recall that in the finite dimensional case
the Clifford multiplication is a deformation of the exterior multiplication
on the exterior algebra.
Therefore such considerations give rise to the possibility of deforming
the cohomology group. 
Vafa suggested to proceed as follows.
The constant maps give an embedding $M \hookrightarrow LM$.
The idea then is to relate the space of differential forms on $M$ 
to a subset of those on $LM$,
take the multiplication there and return to the setting on $M$.

We now give a heuristic justification based on localization 
in equivariant theory,
following the beautiful ideas of Witten in his approach to
index theory \cite{Ati1} and elliptic genera \cite{Wit2}.
Recall that there is a natural action of $S^1$ on $LM$ 
given by rotations on $S^1$
the fixed point set is exactly the set $M$ of constant loops.
Using Fourier series expansion,
one sees that the complexified tangent space of $LM$ restricted to $M$ 
has the following decomposition;
$$TLM|_M \otimes \bC \cong \bigoplus_{n \in \bZ} t^n TM\otimes \bC.$$
The bundle of semi-infinite form on $LM$ restricted to $M$ is
$$\Lambda^{\frac{\infty}{2} + *}(LM)|_M \cong 
\Lambda(\oplus_{n \leq 0} t^n T^*M).$$
When $M$ is endowed with a Riemannian metric,
$\Lambda^{\frac{\infty}{2} + *}(LM)|_M$ is a bundle of conformal vertex 
algebras that contains $\Lambda(T^*M)$ as a subbundle.
This construction can be extended to $LM$.
Now pretend that the exterior differential $d$ has been defined
on $LM$,
we want to consider its localization to $M$.
Witten's work on Morse theory \cite{Wit1} and its generalization by
Bismut \cite{Bis} suggest that $d$ on $LM$ should localize to $d$ on $M$ 
tensored by an algebraic operator on the part 
that comes from the normal bundle.
So on the space of semi-infinite forms restricted to $M$,
we should consider the coupling of $d$ on $M$ with the BRST operator
on the fiber.
Unfortunately the author does not know how to do this.
Nevertheless, 
we can couple the Dolbeault operator with the BRST operator 
when $M$ is complex as 
in \S \ref{sec:complex}, 
\S \ref{sec:Kahler}, 
and \S \ref{sec:hyperkahler}.

So far we have only talked about the fermionic part of the Hilbert space
mentioned in \S 2 of \cite{Vaf}.
To get the bosonic part hence the supersymmetry,
we use the language of supermanifolds (Kostant \cite{Kos}).
Recall that a supermanifold is an ordinary manifold $M$ together with 
a $\bZ_2$-graded structure sheaf.
The even part of the structure sheaf is the sheaf of $C^{\infty}$ function 
on $M$,
while the odd part is the sheaf of sections to the exterior bundle
$\Lambda(E)$ of some vector bundle $E$ on $M$.
The super tangent bundle of $(M, E)$ is a upper vector bundle
$$T(M, E) = TM \oplus E^*,$$
where $TM$ is the even part, $E^*$ is the odd part.
And the differential forms on $(M, E)$ are just sections to 
$\Lambda(T^*M) \otimes S(E)$.
A canonical choice for $E$ is the cotangent bundle $T^*M$,
then we get a supermanifold which corresponds to $\Lambda(T^*M)$.
The supertangent bundle of $(M, T^*M)$ is just
$$T(M, T^*M) = TM \oplus \Pi TM,$$
where $\Pi TM$ means a copy of $TM$ regarded as an odd vector bundle.
We now consider the super loop space $L^sM = Map(S^1, (M, T^*M))$
and regard $(M, T^*M)$ as the fixed point set of the natural circle action.
We have
$$TL^sM|_{(M, T^*M)} \otimes \bC \cong \bigoplus_{n \in \bZ} t^n TM
\bigoplus_{n \in \bZ} t^n \Pi TM.$$
The bundle of semi-infinite form on $L^sM$ restricted to $(M, T^*M)$ is
$$\Lambda^{\frac{\infty}{2} + *}(L^sM)|_{(M, T^*M)} \cong 
\Lambda(\oplus_{n \leq 0} t^n T^*M) \otimes S(\oplus_{n < 0} t^n T^*M).$$
This explains the choice of the vector bundles on $M$ 
for which we construct structures of superconformal vertex algebra bundles.

\medskip

{\bf Acknowledgment}.
{\em The work in this paper is done during the author's visit
at Texas A$\&$M University. 
He appeciates 
the hospitality and the financial support of the Department of Mathematics.
He thanks Catherine Yan for helps on vertex algebras.
Special thanks are due to Huai-Dong Cao and Weiqiang Wang 
for their useful comments on an earlier version.}

\end{document}